\documentclass{conm-p-l}

\newtheorem{theorem}{Theorem}[section]
\newtheorem{lem}[theorem]{Lemma}
\newtheorem{prop}[theorem]{Proposition}
\newtheorem{cor}[theorem]{Corollary}

\theoremstyle{definition}
\newtheorem{defin}[theorem]{Definition}

\theoremstyle{remark}
\newtheorem{rem}[theorem]{Remark}
\newtheorem{ack}{Acknowledgments} 

\numberwithin{equation}{section}

\newcommand{\bbN}{{\mathbb{N}}}
\newcommand{\bbR}{{\mathbb{R}}}
\newcommand{\bbC}{{\mathbb{C}}}
\newcommand{\bbZ}{{\mathbb{Z}}}

\newcommand{\al}{\alpha}
\newcommand{\be}{\beta}

\newcommand{\e}{\varepsilon}

\renewcommand{\span}{\operatorname{span}}
\newcommand{\supp}{\operatorname{supp}}

\newcommand{\sgn}{\operatorname{sgn}}

\newcommand{\lb}{\label}

\newcommand{\lra}{\longrightarrow}

\newcommand{\buo}{without loss of generality }

\newcommand{\INTO}{\buildrel {\mbox{\small into}}\over \longrightarrow}
\makeatletter
\def\@currentlabel{2.1}\label{e:dispaa}
\def\@currentlabel{2.21}\label{e:dispau}
\def\@currentlabel{2.22}\label{e:dispav}
 \def\@currentlabel{2.23}\label{e:dispaw}
\def\@currentlabel{2.24}\label{e:dispax}
\makeatother

\makeatletter
\def\alphenumi{%
  \def\theenumi{\alph{enumi}}%
  \def\p@enumi{\theenumi}%
  \def\labelenumi{(\@alph\c@enumi)}}
\makeatother

\ifx\fiverm\undefined
	\newfont\fiverm{cmr5}
\fi
\input prepictex
\input pictex
\input postpictex
\begin{document}

\title{Injective isometries in Orlicz spaces}

\author{Beata Randrianantoanina}

\address{Department of Mathematics and Statistics
\\ Miami University \\Oxford, OH 45056}

 \email{randrib@@muohio.edu}

%\centerline{\today}

\subjclass{46E30,46B04}

\begin{abstract}
We show that injective isometries in Orlicz space   $L_M$ have to preserve
disjointness, provided that Orlicz function $M$ satisfies
$\Delta_2$-condition, has a continuous second derivative $M''$, satisfies
another ``smoothness type'' condition and either
  $\lim_{t\to0} M''(t) = \infty$ or $M''(0) = 0$ and $M''(t)>0$ for all $t>0$.  
The fact that surjective
isometries of any rearrangement-invariant function space have to preserve
disjointness has been determined before.  However dropping the assumption
of surjectivity invalidates the general method.  In this paper we use a
differential technique.
\end{abstract}
\maketitle

\section{Introduction}

The study of isometries of Banach spaces goes back to Banach's 1932
treatise on linear operators \cite{Ban} and since then it received much
attention in the literature, see the survey \cite{FJs} with its over
300 references.  We just mention here the results most closely related to
the present work.

Banach showed that in separable $L_p(\Omega, \mu), p \neq
2$, any isometry $T: L_p(\Omega, \mu)\lra L_p(\Omega,\mu)$ is of the
form
$$Tf(\omega) = h(\omega)f(\sigma (\omega))$$
where $\sigma: \Omega \lra \Omega$ is a Borel automorphism of $\Omega$ and $h$ is
a scalar function on $\Omega$.  Banach obtained this result by showing that
every isometry $T$ of $L_p$ has to preserve disjointness i.e. if $f,g \in
L_p$ are such that $\mu(\supp f \cap \supp g) = 0$ then also $\mu (\supp Tf
\cap \supp Tg) = 0$ (cf. also \cite{Lam}).
To achieve this he characterized disjointness of functions $f,g$ through
the differential properties of the function $N(\alpha) = \Vert f + \alpha g
\Vert$.  Similar technique was later applied by Koldobsky \cite{Kol91} to
study injective isometries of $L_p(L_q)$ and  Kami\'nska \cite{Kam91},
who observed that isometries of Orlicz
spaces $L_M$ preserve disjointness under assumptions that 
 both
$M$ and $M'$ are strictly convex, $M'(0) = M''(0) = 0$ and $M$ satisfies
$\Delta_2$-condition.

In the present paper we adapt Banach's differential technique to 
Orlicz spaces $L_M$, where $M$ satisfies $\Delta_2$-condition,
another ``smoothness type'' condition  (see Definition~\ref{defdelta2+})
and has a
continuous second derivative $M''$.  We note that unlike Banach and
Koldobsky we do not obtain conditions which are equivalent to the
disjointness of supports of functions $f,g$.  In the case when $M''(0) = 0$
we only describe some conditions which are necessary for disjointness and
some conditions which are sufficient for disjointness, which together
enable us to characterize the injective isometries.  The argument in the
case $\lim_{t \to 0}M''(t) = \infty$ is even more delicate and depends on
the sum of sufficient conditions and necessary conditions for containment
of supports of $f$ and $g$.  The lack of isometric conditions equivalent to
disjointness should not be really surprising in view of the classical
Bohnenblust's characterization of $L_p$-spaces \cite[Theorem~1.b.7]{LT2}.

We conclude the introduction by recalling that surjective isometries of
complex Orlicz spaces were described by Lumer \cite{L62,L63} (reflexive
case) and Zaidenberg \cite{Z76}.  Jamison, Kami\'nska and P.K. Lin \cite{JKL}
 studied surjective isometries of complex
Musielak-Orlicz spaces and real Nakano spaces.  The form of surjective
isometries of real Orlicz spaces follows from the description of surjective
isometries of real rearrangement-invariant spaces \cite{KR}.  
%To the best
%of authors knowledge there has not been any previous work on injective
%isometries of Orlicz spaces other than the very special case in \cite{Kam91}.  

The results of present paper are valid in both complex and real case.

\begin{ack}
I wish to thank Professor A. Koldobsky for introducing me to this problem
and for suggesting the use of the differential technique.  I am also
grateful  to Professors Y. Abramovich and A. Kami\'nska for valuable
discussions, and to Professor D. Burke for help with drawing Figures~1 and 2.

Part of this work was completed while the author participated in the
Workshop in Linear Analysis and Probability at Texas A\&M University,
College Station, Texas in the Summer 1997 organized by Professors W.
Johnson, D. Larson, G. Pisier and J. Zinn.  I wish to thank the organizers
for their hospitality and support.
\end{ack}

\section{Definitions and preliminary lemmas}

%We start with a  generalization of   \cite[Lemma~1]{Kol91}.

We follow standard definitions and notations as may be found e.g. in
\cite{Kr-Rut} or \cite{Chen}.  We recall the basic definitions below.

We say that a function $M: [0, \infty) \lra [0, \infty)$ is an {\it Orlicz
function} if $M$ is convex, $M(0)=0$, $M(1)=1$, $\lim_{u\to0}M(u)/u=0$
and $\lim_{u\to\infty}M(u)/u=\infty$.

The Orlicz function $M$ generates the {\it Luxemburg
norm} defined for scalar valued functions on $\Omega$ by:
$$
\Vert f \Vert_M = \inf\{\lambda: \int_\Omega M \left(\frac{\vert f(u)
\vert}{\lambda}\right)d\mu \leq 1\}.
$$

The {\it Orlicz space} $L_M$ is the space of (equivalence classes of)
measurable functions $f$ with $\Vert f \Vert_{M}<\infty$.

We say that two Orlicz functions
$M_1$ and $M_2$ are {\it equivalent} if there
exist $u_0> 0, k, l> 0$ such that for all
$u> u_0$
$$ M_2(ku)\le M_1(u)\le M(lu). $$

This condition is of importance since
Orlicz spaces $L_{M_1}, L_{M_2}$ are isomorphic
if and only if the Orlicz functions $M_1$, $M_2$ are
equivalent.

It is well known (see e.g. \cite{Chen})
that any Orlicz function $M$ can be ``smoothed out'', that is 
for any $M$ there exists an equivalent
Orlicz function $M_1 $ such that $M_1$ is
twice differentiable and $M_1'' (u)> 0$ for all
$u>0$. 
Moreover, given any $\e> 0$ it is possible
to choose $M_1$ so that $L_M$ and $L_{M_1}$
are $(1+\e)-$isomorphic to each other
\cite{Chen}.

We say that the Orlicz function $M(u)$ {\it satisfies the 
$\Delta_2$ condition
for large values of $u$} if there exist constants $k>0$ and $u_0 \geq 0$
such that for all $u \geq u_0$
$$
M(2u) \leq k M(u).
$$

If the Orlicz function $M$ satisfies the $\Delta_2$ condition then every
 Orlicz function $M_1$ equivalent to $M$ also satisfies the 
$\Delta_2$ condition.

Note that Orlicz function $M$ is convex so it has the right derivative
$M'$.  Krasnoselskii and Rutickii provide the following characterization
of the $\Delta_2$-condition in terms of $M'$.

\begin{prop}
\lb{KRd2}\cite[Theorem 4.1]{Kr-Rut}
A necessary and sufficient condition that the Orlicz function $M(u)$
satisfy the $\Delta_2$-condition is that there exist constants $\alpha$ and
$u_0 \geq0$ such that, for $u \geq u_0$
\begin{equation}\lb{delta2}
\frac{u M'(u)}{M(u)} < \alpha,
\end{equation}
where $M'$ denotes the right derivative of $M$.

Moreover, if \eqref{delta2} is satisfied then $M(2u)\le2^\al M(u)$ for
$u \geq u_0$. 
\end{prop}

Functions which satisfy the $\Delta_2$-condition do not increase more
rapidly than polynomials \cite[p.~24]{Kr-Rut}. In fact there exists a 
constant $C$ so that for all $u\ge u_0$
\begin{equation}\lb{poly}
{M(u)} < Cu^\alpha,
\end{equation}
where $\al$ and $u_0$ are the same as in \eqref{delta2}.

It is also true (see \cite{Chen} or \cite{Kr-Rut}) that for any 
Orlicz function $M$ and for any $u>0$
\begin{equation}\lb{all}
\frac{u M'(u)}{M(u)} >1,
\end{equation}

We now introduce another condition which on
one hand is very similar to \eqref{delta2},
but on the other hand is in  its nature
of ``smoothness type'', as we
explain below.

\begin{defin}\lb{defdelta2+}
Assume that the Orlicz function $M $ is
twice differentiable and that $M $ satisfies  
the $\Delta_2-$ condition.
We say that $M$ {\it satisfies condition
$\Delta_{2+}$} if there exist constants $\be> 0 $ and $u_0\ge 0$
 such that
for all $ u\ge u_0$
\begin{equation}\lb{delta2+}
\frac{u M''(u)}{M'(u)}<\beta.
\end{equation}
\end{defin}

Condition $\Delta_{2+}$ is very important for us because we can prove
our isometry results only for Orlicz spaces $L_M$, such that $M$ satisfies
$\Delta_{2+}$. Thus we wish to observe that this
condition
is of ``smoothness type'' in the following sense: 
\begin{itemize}
\item[(i)] for every function $M$ which satisfies condition $\Delta_{2}$ there
exists an equivalent  Orlicz function
$M_1$ which does  satisfy $\Delta_{2+}$, see Lemma~\ref{lemdelta2+};
However, we do not know whether for every $\e>0$ it is possible to 
choose
$M_1$ so that it is  $(1+\e)-$equivalent with $M$, 
\item[(ii)] for every Orlicz function $M$
which satisfies $\Delta_{2+}$ there
exists an equivalent (even up to an arbitrary $\e>0$) Orlicz function
$M_1$ which does not satisfy $\Delta_{2+}$ (see Lemma~\ref{perturb}).
\end{itemize}

\begin{lem}\lb{lemdelta2+}
Let $M$ be any Orlicz function which
satisfies  condition $\Delta_2$.
Then there exists  an Orlicz function
$M_1$ which is equivalent to $M$
and satisfies condition $\Delta_{2+}$.
\end{lem}

\begin{rem} We do not know whether for every $\e>0$ and every $M$
satisfying $\Delta_2$ it is possible to choose $M_1$ 
 satisfying  $\Delta_{2+}$ and $(1+\e)-$equivalent to $M$.
\end{rem}

\begin{proof}
Define $M_1'$ by (see Figure~1):
\begin{equation*}
M_1'(u)=
\begin{cases}
M'(u) \hspace{45mm} \text{if} \  u\le2\\
f_k(u) \hspace{46mm} \text{if} \  2^k\le u\le 2^k+1\ ,  k\in\bbN\\
M'(2^k)+(u-2^k)\frac{M'(2^{k+1})-M'(2^k)}{2^k}
 \hspace{3mm} \text{if} \  2^k+1<u<2^{k+1} \ ,  k\in\bbN
\end{cases}
\end{equation*}
where $f_k:\bbR\lra\bbR$ is an increasing, continuous, differentiable
function such that $f_k'$ is monotone on $[2^k,2^k+1]$ and
\begin{eqnarray*}
f_k(2^k)&=&M'(2^k)\\
f_k(2^k+1)&=&M'(2^k)+\frac{M'(2^{k+1})-M'(2^k)}{2^k}\\
f_k'(2^k)&=&\frac{M'(2^{k})-M'(2^{k-1})}{2^{k-1}}\\
f_k'(2^k+1)&=&\frac{M'(2^{k+1})-M'(2^k)}{2^k}.
\end{eqnarray*}
%\input graph4
%%Version: Dec 1,1998 11:12:22
%%Graph of y=2.1*x/3 +sin(x)/2+sin(x)*cos(1+2*x)/4 from 0 to 10
%\input prepictex
%\input pictex
%\input postpictex
%\input latexpicobjs %See docs
%\magnification=1200

\beginpicture
%\startrotation by 0 -1 about 0 0 %%% -90 degrees
\setcoordinatesystem units <1cm,.7cm> point at 0 10
%%Adjust next line for circles.
%\circulararc 360 degrees from 2 0 center at 0 0

\setquadratic
\plot
0.0000	 0.0000
0.2500	 0.3031
0.5000	 0.5398
0.7500	 0.7293
1.0000	 0.9125
1.2500	 1.1273
1.5000	 1.3857
1.7500	 1.6651
2.0000	 1.9191
2.2500	 2.1019
2.5000	 2.1929
2.7500	 2.2090
3.0000	 2.1972
3.2500	 2.2115
3.5000	 2.2874
3.7500	 2.4252
4.0000	 2.5940
4.2500	 2.7506
4.5000	 2.8663
4.7500	 2.9442
5.0000	 3.0195
5.2500	 3.1418
5.5000	 3.3484
5.7500	 3.6441
6.0000	 3.9969
6.2500	 4.3535
6.5000	 4.6649
6.7500	 4.9101
7.0000	 5.1037
7.2500	 5.2852
7.5000	 5.4944
7.7500	 5.7476
8.0000	 6.0266
8.2500	 6.2869
8.5000	 6.4811
8.7500	 6.5841
9.0000	 6.6079
9.2500	 6.5965
9.5000	 6.6048
9.7500	 6.6716
10.0000	 6.8025
/
%10.25 7/
%\setplotarea x from 0 to 10, y from -2 to 2
%\put {$y=2.1*x/3 +sin(x)/2+sin(x)*cos(1+2*x)/4$} [r] at 10 2
\put {$\bullet$} at 1.2500 1.1273
\put {$\bullet$} at 2.5000 2.1929
\put {$\bullet$} at 5.0000 3.0195
\put {$\bullet$} at 10.0000 6.8025
\put {1} <0mm, -4mm> at .625 0
\put {2} <0mm, -4mm> at 1.2500 0
\put {$2^2$} <0mm, -4mm> at 2.500 0
\put {$2^3$} <0mm, -4mm> at 5.00 0
\put {$2^4$} <0mm, -4mm> at 10.00 0
\put {$M'$} <-2mm, 2mm> at 8.2500 6.2869
\put {$M_1'$} <2mm, -2mm> at 7.0000 4.5
\put {Figure~1.}  at 5.0000 -2
\put { }  at -1.000 7.5
\put { }  at -1.000 -2.5
\setlinear
\plot 0 0 10.5 0 / %% x-axis
\plot 0 -1 0 7 / %% y-axis
\plot 1.2500 1.1273 2.5000 2.1929 /
\plot 2.5000 2.1929 5.0000 3.0195 /
\plot 5.0000 3.0195 10.0000 6.8025 /
%%Some ticks:
\plot .625 .03 .625 -.03 / %Along the x-axis
\plot 1.25 .03 1.25 -.03 /
\plot 2.5 .03 2.5 -.03 /
\plot 5 .03 5 -.03 /
\plot 10 .03 10 -.03 /
%%\plot 2 .03 2 -.03 /
%%\plot 3 .03 3 -.03 /
%%\plot 4 .03 4 -.03 /
%\plot .03 2 -.03 2 / %Along the y-axis
%\plot .03 1 -.03 1 /
%\plot .03 -2 -.03 -2 /
%\plot .03 -1 -.03 -1 /

%%\stoprotation
\endpicture

Thus $M_1'(u)$ is  increasing, continuous and differentiable for all
$u\in\bbR_+$.
Then 
$$M_1(u)=\int_0^u M_1'(t)\ dt.$$

We first check that $M_1$ is equivalent to $M$. Indeed,
for any $u\ge2$ let $k\in\bbN$ be such that $2^k\le u<2^{k+1}$.
Then we have:
\begin{equation*}
\begin{split}
M_1(u)&=\int_0^u M_1'(t)\ dt
\ge \int_0^2 M'(t)\ dt + \sum_{j=1}^k 2^jM'(2^j)\\
&\text{by \eqref{all}}\\
&>\sum_{j=1}^k M(2^j)>M(2^k)\ge M(\frac u2)
\end{split}
\end{equation*}

Also we have:
\begin{equation*}
\begin{split}
M_1(u)&=\int_0^u M_1'(t)\ dt\\
&\le \int_0^2 M'(t)\ dt + \sum_{j=1}^k 2^jM'(2^{j+1})+(u-2^k)M'(2^{k+1})\\
&\text{by \eqref{delta2}}\\
&\le M(1)+\frac\al2\sum_{j=1}^{k+1} M(2^{j+1})
\le \frac12+\frac\al2\sum_{l=1}^{k+2} M(2^{l})\\
&\text{by convexity of $M$}\\
&\le \frac12+\frac\al2\sum_{l=1}^{k+2} \frac{M(2^{k+2})}{2^{k+2-l}}\\
&\le \frac12+\frac\al2 M(2^{k+2})
\le \frac12+\frac\al2 M(4u)<\al M(4u)\\
&\text{since $\al>1$}\\
&\le M(4\al u)
\end{split}
\end{equation*}

Since for $u\le2$ \ $M(u)=M_1(u)$, we conclude that for all $u\ge0$
$$M(\frac u2)\le M_1(u)\le M(4\al u)$$
i.e. $M_1$ is equivalent to $M$.

Next we show that $M_1$ satisfies condition $\Delta_{2+}$.

Let $u\ge 2$ and $k\in\bbN$ be such that $2^k\le u<2^{k+1}$.
Then, if $2^k+1\le u$ we have:
\begin{equation}\lb{main}
\begin{split}
\frac{uM_1''(u)}{M_1'(u)}&\le\frac{2^{k+1}\frac{M'(2^{k+1})-M'(2^k)}{2^k}}
{M'(2^k)}=\frac{2M'(2^{k+1})}{M'(2^k)}-2\\
&\text{ by \eqref{delta2} and \eqref{all}}\\
&\le\frac{2\frac{\al M(2^{k+1})}{2^{k+1}}}{\frac{M(2^{k})}{2^{k}}}
=\frac{\al M(2^{k+1})}{M(2^{k})}\\
&\text{ by Proposition~\ref{KRd2}}\\
&\le \al2^\al.
\end{split}
\end{equation}

If $2^k\le u\le2^k+1$, then $M_1'(u)=f_k(u)$ and we consider
two cases:

$(1)$ If $f_k'$ is increasing on $[2^k,2^k+1]$ then
$$f_k'(u)\le f_k'(2^k+1)=\frac{M'(2^{k+1})-M'(2^k)}{2^k}$$
and
\begin{equation*}
\begin{split}
\frac{uM_1''(u)}{M_1'(u)}&=\frac{uf_k'(u)}{f_k(u)}\le
\frac{2^{k}+1\frac{M'(2^{k+1})-M'(2^k)}{2^k}}{M'(2^k)}\\
&\text{ by \eqref{main}}\\
&\le \al2^\al.
\end{split}
\end{equation*}

$(2)$ If $f_k'$ is decreasing on $[2^k,2^k+1]$ then
$$f_k'(u)\le f_k'(2^k)=\frac{M'(2^{k})-M'(2^{k-1})}{2^{k-1}}$$
and
\begin{equation*}
\begin{split}
\frac{uM_1''(u)}{M_1'(u)}&=\frac{uf_k'(u)}{f_k(u)}\le
\frac{(2^k+1)\frac{M'(2^{k})-M'(2^{k-1})}{2^{k-1}}}{M'(2^k)}\\
&\le
\frac{2\cdot 2^{k}\frac{M'(2^{k})-M'(2^{k-1})}{2^{k-1}}}{M'(2^{k-1})}\\
&\text{ by \eqref{main} applied to $(k-1)$}\\
&\le 2\al2^\al.
\end{split}
\end{equation*}
\end{proof}

We now turn our attention to property
(ii) of $\Delta_{2+}$. Note that (ii) is trivially true because
of the requirement that every function
satisfying $\Delta_{2+}$ is twice
differentiable, which can be easily
perturbed. But even more is true.

\begin{lem}\lb{perturb}
Let $M$ be an Orlicz function which satisfies
condition $\Delta_{2+}$ and let 
$e> 0$ be  given. Then there exists
  a twice
differentiable Orlicz function $M_1$ such that
$$ 
(1-\e)M(u) \le M_1(u)\le (1+\e)M(u)$$
for all $u\in\bbR$ and which does not satisfy condition
$\Delta_{2+}$.
\end{lem}

{\sc Sketch of proof.}
This fact is not of a particular
importance for our paper 
%(in fact we would prefer that it wasn't true) 
so we do
not provide a detailed proof but only
sketch the idea.

Let $M$ be an Orlicz function which satisfies \eqref{delta2}
and \eqref{delta2+}.
For any $n\in\bbN$ put
\begin{equation*}
\begin{split}
\al(n) &= M'(n+\frac12+\frac1{2^{2n}})-M'(n+\frac12-\frac1{2^{2n}})\\
\be(n) &=\frac{\al(n)}{2^n}
\frac{n+\frac12-\frac1{2^{2n}}}{M'(n+\frac12+\frac1{2^{2n}})}.
\end{split}
\end{equation*}

Then there exists $n_0\in\bbN$ so that for all $n\ge n_0$
$$\be(n)<\frac12\frac1{2^n}.$$

Indeed, by Mean Value Theorem, for each $n$ there exists
$c(n)\in(n+\frac12-\frac1{2^{2n}},n+\frac12-\frac1{2^{2n}})$
so that
$M''(c(n))=\al(n)\cdot2^n.$ Then, by \eqref{delta2+},
$$
\be >\frac{c(n)M''(c(n))}{M'(c(n))}=2^n\al(n)\frac{c(n)}{M'(c(n))} \ge\be(n)\cdot2^{2n}.
$$

Also this calculation implies that
\begin{equation*}
\begin{split}
\al(n)&<\frac{\be M'(c(n))}{2^nc(n)}\\
&{\text{ by \eqref{delta2}}}\\
&<\frac{\be\al M(c(n))}{2^n(c(n))^2}\\
&{\text{ by \eqref{poly}}}\\
&\le \frac{\be\al C(c(n))^\al}{2^n(c(n))^2}\\
&\le \frac{\be\al C(n+1)^\al}{2^n n^2}.
\end{split}
\end{equation*}

Thus $\al(n)\lra 0$  as $n\lra\infty$ and, for a given
$\e>0$ we can choose $n_1\ge n_0$ so that for all $n\ge n_1$
$$\al(n)<\e.$$
For $n\ge n_1$ we ``adjust'' $M'$ on $(n,n+1)$ by setting
(see the graph in
Figure~2):
\begin{equation*}
\tilde{M}'(x)=
\begin{cases}
M'(n+\frac12-\frac1{2^{2n}}) \hspace{2mm}\text{ if } 
x\in(n+\frac12-\frac{1}{2^{2n}},
n+\frac12-\frac{\be(n)}{2}) \text{ for some } n\ge n_1,\\
M'(n+\frac12+\frac{1}{2^{2n}}) \hspace{2mm}\text{ if } 
x\in(n+\frac12+\frac{\be(n)}2,n+\frac12+\frac1{2^{2n}}) 
\text{ for some } n\ge n_1,\\
M'(n+\frac12-\frac1{2^{2n}}) +(x-(n+\frac12-\frac{\be(n)}2))\cdot
\frac{\al(n)}{\be(n)} \\
\hspace{27mm}\text{ if } x\in(n+\frac12-\frac{\be(n)}{2},
n+\frac12+\frac{\be(n)}2) \text{ for some } n\ge n_1,\\
M'(x) \hspace{18mm}\text{ otherwise.}
\end{cases}
\end{equation*}

%\input graph3
%%Version: Dec 1,1998 11:09:32
%%Graph of y=2.3*x/3 +sin(.5*x)/3+sin(x)/3+cos(1+2*x)/5 from 0 to 10
%\input prepictex
%\input pictex
%\input postpict
%\input latexpicobjs %See docs
%\magnification=1200
%\input pictex.sty

\beginpicture
%\startrotation by 0 -1 about 0 0 %%% -90 degrees
\setcoordinatesystem units <.8cm,.8cm> point at -5 5
%%Adjust next line for circles.
%\circulararc 360 degrees from 2 0 center at 0 0

\setquadratic
\plot
0.0000	 0.1081
0.2500	 0.3298
0.5000	 0.5424
0.7500	 0.7641
1.0000	 1.0090
1.2500	 1.2824
1.5000	 1.5790
1.7500	 1.8834
2.0000	 2.1737
2.2500	 2.4268
2.5000	 2.6245
2.7500	 2.7578
3.0000	 2.8303
3.2500	 2.8578
3.5000	 2.8653
3.7500	 2.8821
4.0000	 2.9353
4.2500	 3.0440
4.5000	 3.2157
4.7500	 3.4447
5.0000	 3.7141
5.2500	 4.0000
5.5000	 4.2775
5.7500	 4.5263
6.0000	 4.7354
6.2500	 4.9051
6.5000	 5.0463
6.7500	 5.1769
7.0000	 5.3168
/

\put{$\updownarrow$} at  1.1 2.9201
\put{$\varepsilon>\alpha(n)$} at  -.1 2.9201
\put{$\leftrightarrow$} at  3.5 1.05
\put{$\beta(n)$} at  3.5 .7
\put{$\leftarrow$} at  2.7 -.6
\put{$\rightarrow$} at  4.2 -.6
\put{${\displaystyle {\frac{1}{2^n}}}$} at  3.5 -.6
\put{${\displaystyle {n+{\frac1{2}}}}$} at  3.5 -1.8
\put{$n$} at   .5 -1.6
\put{$n+1$} at   6.5 -1.6
\put{$M'$} at   6.7 4.8
\put{$\searrow$} at  3.7 3.5
\put{$\tilde{M}'$} at   3.4 3.9
\put{Figure~2.} at   3.5 -3
\put{$ $} at   -5 -3.5
%\setplotarea x from 0 to 10, y from -2 to 2
%\put {$y=2.3*x/3 +sin(.5*x)/3+sin(x)/3+cos(1+2*x)/5$} [r] at 10 2
\setlinear
\plot 0 -1.2 7 -1.2 / %% x-axis
\plot 2.5000 2.6245 3.25 2.6245 /
\plot 3.25 2.6245 3.75 3.2157 /
\plot 3.75 3.2157 4.5000 3.2157 /
%\plot 0 -2 0 2 / %% y-axis
%%Some ticks:
\plot .5 -1.23 .5 -1.17 /
\plot 6.5 -1.23 6.5 -1.17 /
\plot 3.5 -1.23 3.5 -1.17 /
\setdots
\plot 3.75 1.0 3.75 3.2157  / %dotted
\plot 4.5 -1 4.5000 3.2157 / %dotted
\plot 1 2.6245 2.5000 2.6245 / %dotted
\plot 1.0 3.2157 3.75 3.2157  / %dotted
\plot 2.5 -1 2.5000 2.6245 / %dotted
\plot 3.25 1.0 3.25 2.6245 / %dotted
%\setcoordinatesystem units <1cm,1cm>
%\betweenarrows {p}  <3mm,-5mm> from 13.75 -1 to 14.25 -1
%\put(.45,-1.3){n} /
%\plot -1 .03 -1 -.03 / %Along the x-axis
%\plot 1 .03 1 -.03 /
%%\plot 2 .03 2 -.03 /
%%\plot 3 .03 3 -.03 /
%%\plot 4 .03 4 -.03 /
%\plot .03 2 -.03 2 / %Along the y-axis
%\plot .03 1 -.03 1 /
%\plot .03 -2 -.03 -2 /
%\plot .03 -1 -.03 -1 /

%%\stoprotation
\endpicture

$\tilde{M}'$ is not differentiable at countable  number of points
(four on each interval $(n,n+1)$ when $n\ge n_1$), but it can be further
adjusted to produce $M_1'$ so that $M_1'$ is differentiable
everywhere, $|M_1'(x)-M'(x)|<\e$ for all $x\ge0$ and 
$$M_1'(n+\frac12)=\tilde{M}'(n+\frac12),$$
$$M_1''(n+\frac12)=\tilde{M}''(n+\frac12)=\frac{\al(n)}{\be(n)}.$$

We leave up to an interested reader the detailed formula for
$M_1$ and the check that $M_1$ is $(1+\e)-$equivalent with $M$.

For $M_1$ we have for $n\ge n_1$
$$\frac{(n+\frac12)M_1''(n+\frac12)}{M_1'(n+\frac12)}=
\frac{(n+\frac12)\al(n)M'(n+\frac12+\frac1{2^{2n}})}
{M_1'(n+\frac12)\al(n)(n+\frac12-\frac1{2^{2n}})}\ge2^n.$$
Hence $M_1$ does not satisfy condition $\Delta_{2+}$.\hfill\qed

In the case when $M$ satisfies the $\Delta_2$-condition and $\Vert f
\Vert_M < \infty$ we have
$$
\int_\Omega M \left(\frac{\vert f(u) \vert}{\Vert f \Vert_{M}}\right) d\mu = 1.
$$

When $M$ satisfies condition $\Delta_2$ then the dual space of $L_M$ is
also an Orlicz space, which is determined by an Orlicz function $M^*$
called {\it a complementary Orlicz function to $M$} which is defined by
$$ M^*(v)=\int_0^v q(s)\ ds$$
for $v\ge 0$, where $q$ denotes the right inverse of $M'$ -- the right
derivative of $M$.

It is well-known (see \cite{Kr-Rut}) that $M^*$ does not have to satisfy
condition $\Delta_2$ whenever $M$ does  and also vice-versa: $M^*$ may 
satisfy
condition $\Delta_2$ when $M$ does not.
Notice, however, 
that when $M$ is twice differentiable, $M''$ is continuous and $M''(t)>0$
for all $t>0$ then the same is true for $M^*$, i.e. $M^*$ is 
twice differentiable, ${M^*}''$ is continuous and ${M^*}''(t)>0$
for all $t>0$. Moreover, in this situation, if $\lim_{t\to 0}M''(t)=\infty$
then $\lim_{t\to 0}{M^*}''(t)=0$.

For any $f\in L_M$ we will denote by $f^N$ 
the {\it norming
functional} for $f$ i.e. the functional  such that
$\|f\|_{L_M}=\|f^N\|_{L_M^*}$ and $f^N(f)=\|f\|_{L_M}^2$.

We will need the following simple fact:

\begin{lem}\lb{calvert}
Let $M$ be a differentiable Orlicz function which satisfies
condition $\Delta_2$, $A, B$ be disjoint subsets of
$[0,1]$ and $H = \span \{\chi_A, \chi_B\}$.
Suppose that $T:L_M\INTO L_M$ is an isometry.

Then $(TH)^* = \{(Th)^N\ :\ h\in H\}=\span \{(T\chi_A)^N, (T\chi_B)^N\}$
and the map $S:\span \{\chi_A^N,
\chi_B^N\}\lra\span \{(T\chi_A)^N, (T\chi_B)^N\}$ defined by
$$S\chi_A^N=(T\chi_A)^N,\ \ \ \ \  S\chi_B ^N=(T\chi_B)^N$$
is an isometry.
\end{lem}

\begin{proof}
Since $M$ is differentiable and satisfies
condition $\Delta_2$, by \cite{GrzH}, $L_M$ is smooth.
Since $L_M$ is a symmetric function space,
for every $A\subset[0,1]$, the norming functional
$\chi_A^N = a \chi_A$ for an appropriate constant $a$.
Indeed, if $\sigma$ is a homeomorphism of $[0,1]$ with
$\sigma(A)= A$, then $Tf=f\circ \sigma$ is a surjective
isometry of $L_M$ such that $T\chi_A=\chi_A$ and
$(T^{-1})^*\chi_A^N=(T\chi_A)^N=\chi_A$. Thus 
$$\chi_A^N=\chi_A^N\circ\sigma^{-1}$$
for any homeomorphism of $[0,1]$ with
$\sigma(A)= A$. Hence $\chi_A^N\big|_A= const.$
 
Similarly, for any disjoint subsets $A,B\subset[0,1]$,
by considering homeomorphisms of $[0,1]$ with
$\sigma(A)= A$ and $\sigma(B)= B$ and the isometries of $L_M$
that they induce, we conclude that for every 
$h \in H=\span\{\chi_A,\chi_B\}$, 
the norming
functional for $h$  
is a linear combination of
norming functionals $\chi_A^N, \chi_B^N$ and $H^* = \span \{\chi_A^N,
\chi_B^N\}$.

Let $V = T^{-1}: TL_M \lra L_M$.
Since $L_M$ is smooth  for all $f\in L_M$ we have:
$$(Tf)^{N} = V^*(f^N).$$
Thus for every $h \in H$:
\begin{equation*}
(Th)^{N} = V^*(h^{N})\in \span \{V^*(\chi_A^N), V^*(\chi_B^N)\}= 
 \span \{(T\chi_A)^N, (T\chi_B)^N\}
\end{equation*}
and
$$(TH)^* = \span \{(T\chi_A)^N, (T\chi_B)^N\}.$$

Further $V^*: H^* \lra (TH)^*$ is an isometry between subspaces of $L_{M^*}$
with
$$V^* \chi_A^N = (T\chi_A)^N \ \ ,\ \ V^*\chi_B^N = (T\chi_B)^N.$$
\end{proof}

We finish this section with a lemma about differentiability of the function
\begin{equation}\lb{F}
F(\alpha,\eta) =
\int M\left(\frac{|f+\alpha g|}{\eta}\right)\,d\mu (t) - 1,
\end{equation}
 where $\al\in\bbR,$ $\eta>0$ and $f,g$ are given functions
of norm $1$ from $L_M$. We will need this lemma to
describe  the differential behaviour of 
$N(\alpha) = \Vert f + \alpha g \Vert_M$.
Lemma~\ref{kol} is inspired by and generalizes
\cite[Lemma~1]{Kol91}.

\begin{lem}\label{kol}
Let $M$ be an Orlicz function satisfying condition $\Delta_{2+}$
and such that $M''$ is continuous. Suppose that
$f,g\in L_M$ with $\|f\|= \|g\|=1$. For $\al\in\bbR,$ $\eta>0$, consider 
the function $F$ defined by \eqref{F}

Then
\begin{itemize}
\item[(a)]
${
\frac{\partial^2 F}{\partial\alpha^2} (\alpha,\eta)}$ is a 
continuous function
with respect to $\eta$, when $\eta\in(0,1)$ and
$$
\frac{\partial^2 F}{\partial\alpha^2} (\alpha,1)
 = \int M''  ({|f+\alpha g| })
{g^2}\,d\mu (t)\ $$
for a.e. $\alpha \in \bbR$,
\item[(b)]
%${\disp
%\frac{\partial F}{\partial \al} (\alpha,\eta)}$ and
${
\frac{\partial F}{\partial \eta} (\alpha,\eta)}$ is a continuous function
with respect to both variables and
%${\disp
%\big|{\partial F\over\partial \al} (0,1)\big|<\infty}$ and
$${
0<\big|\frac{\partial F}{\partial \eta} (0,1)\big|<\infty,}$$
\item[(c)]
${
\frac{\partial^2 F}{\partial\al\partial \eta} (\alpha,\eta)}$ 
and ${
\frac{\partial^2 F}{\partial \eta^2} (\alpha,\eta)}$ are continuous
functions
with respect to both variables and
%${\disp
%\big|\frac{\partial F}{\partial \al} (0,1)\big|<\infty}$ and
$${
\big|\frac{\partial^2 F}{\partial\al\partial \eta} (0,1)\big|<\infty,\ \ \ 
\ \ \ \ \ \ 
\big|\frac{\partial^2 F}{\partial \eta^2} (0,1)\big|<\infty.}$$
\end{itemize}
\end{lem}

\begin{proof}
The proofs of parts $(a), (b), (c)$ are very similar to each other and
essentially consist of an application of Fubini's and Lebesgue's theorems.

We will need the following auxiliary functions,
for  $\alpha \in \bbR$, $\eta>0$:
\begin{eqnarray*}
h_\eta(\alpha) & = &\int M'  \left(\frac{|f+\alpha g|}{\eta}\right)
\sgn (f+\alpha g) \frac{g}{ \eta}\,d\mu (t)\ ,\cr
s_\eta(\alpha) & = & \int M''  \left(\frac{|f+\alpha g|}{ \eta}\right)
\frac{g^2}{\eta^2}\,d\mu (t)\ ,\cr
w_\al(\eta) & = &\int M' \left( \frac{|f+\alpha g|}{\eta}\right)
\left(- \frac{|f+\alpha g|}{\eta^2}\right)\, d\mu (t)\ ,\cr
v(\al,\eta) & = & -\int M''  \left(\frac{|f+\alpha g|}{ \eta}\right)
\frac{(f+\alpha g)g}{\eta^3}\,d\mu (t)\\
& \hspace{3mm} & \hspace{5mm} -\int M' \left( \frac{|f+\alpha g|}{\eta}\right)
\sgn (f+\alpha g) \frac{g}{\eta^2}\,d\mu (t),\\
z_\al(\eta) & = & \int M''  \left(\frac{|f+\alpha g|}{ \eta}\right)
\frac{|f+\alpha g|^2}{\eta^4}\,d\mu (t)\\
& \hspace{3mm} & \hspace{5mm} 
+\int M' \left( \frac{|f+\alpha g|}{\eta}\right)
\frac{2|f+\alpha g|}{\eta^3}\, d\mu (t)\ 
\end{eqnarray*}

First notice that, since $M$ satisfies condition \eqref{delta2+},
there exists $\be>0$ so that for all $\al\in\bbR$
$$M''\left(\frac{|f+\al g|}{\eta}\right)|f+\al g|\le\be
M'\left(\frac{|f+\al g|}{\eta}\right).$$

By  \cite[Lemma~9.1, p.73]{Kr-Rut},
for any $\alpha \in \bbR$ and $\eta>0$
$$M' \left( \frac{|f+\alpha g|}{\eta}\right) \in L_{M^*}\ ,$$
where $M^*$ is the complementary Orlicz function to $M$.
Thus, by the generalized H\"older's inequality
(\cite[Theorem~9.3, p.74]{Kr-Rut}):
\begin{equation} \lb{bdd}
 |h_\eta(\alpha)| <\infty\ ,
 |w_\al(\eta)| < \infty \ , |v(\al,\eta)| < \infty\ ,
 |z_\al(\eta)| < \infty, 
 \end{equation}
  for all
$\alpha\in\bbR$, $\eta>0$.

Notice that for any $\beta, u,v\in\bbR$, $\eta>0$ we have:
$$\int_0^\beta M'' \left(\frac{|u+\alpha v|}{\eta}\right) \frac{v^2}{\eta^2}
\,d\alpha
= M' \left(\frac{|u+\alpha v|}{\eta}\right) \sgn (u+\al v)
\frac{v}{\eta}\big|_0^\beta$$
So, by the Fubini theorem
$$\int_0^\beta s_\eta(\alpha) \,d\alpha
= h_\eta(\beta) - h_\eta(0)\ .$$
Thus $s_\eta$ is  absolutely integrable on $[0,\beta]$ with respect
to $\al$, $|s_\eta|<\infty$ for almost all $\al$ and $h_\eta$ is
a primitive for $s_\eta$.
By Lebesgue's theorem, for each $\eta>0$,
$$\frac{d}{d\al}h_\eta (\alpha) = s_\eta(\alpha)$$
 for almost
all $\alpha\in\bbR$.

Similarly, for each $\eta>0$,
 $\int_0^\beta h_\eta(\alpha)\,d\alpha = F(\alpha,\eta)|_0^\beta$
and
$$\frac{\partial F}{\partial\alpha} (\alpha,\eta) =
h_\eta(\alpha)$$
 for almost
all $\alpha \in \bbR$.
Thus for $\eta=1$
$$s_1(\al) = \frac{d}{d\al}h_1 (\al) =
\frac{\partial}{\partial\alpha}\frac{\partial F}{\partial\alpha} (\alpha,\eta)$$
and   we get the formula in part $(a)$.

The continuity of ${
\frac{\partial^2 F}{\partial\alpha^2} (\alpha,\eta)}$ 
follows from the fact that $M''$ is continuous (see e.g.
\cite[Sections~IV.2, IV.4]{Courant}).

To prove $(b)$ we use the same argument after we notice that
for each
$\alpha \in \bbR$ and for any $\varepsilon,\zeta >0$:
$$\int_\varepsilon^\zeta M' \left( \frac{|u+\alpha v|}{\eta}\right)
(- \frac{|u+\alpha v|}{\eta^2})\, d\eta =
M \left( \frac{|u+\alpha v|}{\eta}\right)\Big|_\varepsilon^\zeta.$$

Thus, as above, $\frac{\partial F}{\partial \eta} (\alpha,\eta) = w_\al(\eta)$
for almost  all $\eta>0$ and the continuity of ${
\frac{\partial F}{\partial\eta} (\alpha,\eta)}$ 
follows from the continuity  of $M'$ \cite{Courant}.

Since $w_\al(\eta) <0$ for all $\eta>0$, we get
$\frac{\partial F}{\partial \eta} (\alpha,1) \ne0$.

Similarly, to get $(c)$ we repeat the argument from above since we have for
each $\eta > 0$, and $\varepsilon, \zeta > 0$
$$\int_\varepsilon^\zeta \left( - M''\left( \frac{|u+\alpha v|}{\eta}\right)
 \frac{(u + \alpha v)v}{\eta^3}-
M' \left( \frac{|u + \alpha v|}{\eta} \right) 
\sgn{(u + \alpha v)}\frac{v}{\eta^2} \right) {d \alpha}$$
$$= \left( - M' \left( \frac{|u + \alpha v|}{\eta} \right) 
\frac{|u + \alpha
v|}{\eta^2} \right) \Big|_{\al=\varepsilon}^{\al=\zeta}$$

Thus 
$$v (\alpha, \eta) = \frac{\partial}{\partial \alpha} w_\alpha (\eta) =
\frac{\partial^2 F}{\partial \alpha \partial\eta} (\alpha, \eta)$$
 for almost all
$\alpha \in \bbR, \eta > 0.$

Finally, for
each $\eta > 0$, and $\varepsilon, \zeta > 0$
$$\int_\varepsilon^\zeta \left(M''\left( \frac{|u+\alpha v|}{\eta}\right)
 \frac{|u + \alpha v|^2}{\eta^4}+
M' \left( \frac{|u + \alpha v|}{\eta} \right) 
\frac{2|u + \alpha v|}{\eta^3}\right) {d \eta}$$
$$= \left( - M' \left( \frac{|u + \alpha v|}{\eta} \right) 
\frac{|u + \alpha
v|}{\eta^2} \right) \Big|_{\eta=\varepsilon}^{\eta=\zeta}$$

Thus 
$$z_\alpha(\eta) = \frac{\partial}{\partial \eta} w_\alpha (\eta) =
\frac{\partial^2 F}{\partial\eta^2} (\alpha, \eta)$$
 for almost all
$\alpha \in \bbR, \eta > 0.$

Continuity of both $v (\alpha, \eta) $ and $z_\alpha(\eta)$ 
with respect to $\alpha$ and $\eta$ is again a consequence of
continuity of $M'$ and $M''$.  The final statement follows from
\eqref{bdd}.
\end{proof}

\section{The case of $M''(0)=0$}

We first study the case of Orlicz spaces $L_M$ analogous to $L_p,p>2$ in a
sense that $M''(0)=0$.  We obtain the following partial description of
functions with disjoint supports.

%This is where Jean began
\begin{prop}\lb{Mzero}
Assume that $M$
is an Orlicz function which satisfies condition
$\Delta_{2+}$  and such that
$M''$ is a continuous function with $M''(0)=0$
and $M''(t) > 0$ for all $t > 0$.
Let $f,g\in L_M$ and $N(\alpha) = \|f+\alpha g\|_M$.

Then
\begin{enumerate}
\item[(a)] If $f, g $
have disjoint supports and $g$ is bounded then $N'(0) = 0$ and
$N''(\alpha) \lra 0$ as $\al\lra 0$ along a subset of $[0,1]$
of full measure.

\item[(b)] If $N'(0) = 0$ and
$N''(\alpha) \lra 0$ as $\al\lra 0$ along a subset of $[0,1]$
of full measure then $f, g $
have disjoint supports.
\end{enumerate}
\end{prop}

\begin{proof}
First note that since $M$ is differentiable and satisfies condition
$\Delta_2$ thus, by \cite{GrzH}, $L_M$ is smooth and the function
$N(\al)$ is differentiable. Since $N(\al)$ is a convex function
of $\al$ also the second derivative $N''(\al)$ exists a.e.

Notice that if $f, g \in L_M$ are disjointly supported then $N$ clearly has
a minimum at $0$, so
$N' (0)=0$.  Thus in the following we will work under the assumption that
$N'(0) = 0$.

Assume, \buo that $\|f\|_M=\|g\|_M=1$, and let $F(\alpha,\eta)$
be defined as in Lemma~\ref{kol}.  Since $M$ satisfies
$\Delta_2$-condition we have
$F(\alpha, N(\alpha)) = 0$ for all $\alpha \in \bbR$.   Therefore
$\frac{d}{d\alpha}(F(\alpha, N(\alpha))) = 0$ for all $\alpha$.  Hence

\begin{equation}\lb{der0}
0 = \frac{d}{d\alpha} (F(\alpha, N(\alpha))) = 
\frac{\partial F}{\partial\alpha}(\alpha,N(\alpha)) + 
\frac{\partial F }{\partial\eta}(\alpha,N(\alpha)) N'(\alpha)
\end{equation}

and by taking the derivative again we get:

\begin{equation}\lb{der}
\begin{split}
- \frac{\partial^2F}{\partial\alpha^2} (\alpha,N(\alpha)) &=
\frac{\partial^2F}{\partial\alpha\partial\eta}
(\alpha,N(\alpha))N'(\alpha) \\
&+\left[\frac{\partial^2F}{\partial\alpha\partial\eta}
(\alpha,N(\alpha))+\frac{\partial^2F}{\partial\eta^2}
(\alpha,N(\alpha))N'(\alpha)\right]N'(\alpha)\\
&+ \frac{\partial F}{\partial\eta}(\alpha,N(\alpha))N''(\alpha)
\end{split}
\end{equation}

Let $A$ be the set so that for $\alpha \in A$,  $N'(\al)$ exists and
$\frac{\partial^2F}{\partial\alpha^2}(\alpha,1)$
is given by the formula from Lemma~\ref{kol}(a).  By Lemma~\ref{kol}(c)
there exists a sequence
$(\alpha_n) \subset A, \alpha_n \to 0$ so that
$\frac{\partial^2F}{\partial\alpha \partial\eta}(\alpha_n, N(\alpha_n))$ 
and 
$\frac{\partial^2F}{\partial\eta^2}(\alpha_n, N(\alpha_n))$ are
bounded for all $n \in \bbN$. Since $N'(\alpha_n) \to 0$ and since, by
Lemma~\ref{kol}(b), ${ \frac{\partial F}{\partial \eta}
(\alpha_n,N(\alpha_n))
\not= 0}$ for $n$ large enough, we get that

\begin{equation*}
\lim_{n\to\infty} \frac{\partial^2F}{\partial\alpha^2}(\alpha_n, N(\alpha_n))
= 0 \Longleftrightarrow \lim_{n \to \infty} N''(\alpha_n) = 0 \ .
\end{equation*}

Since $\frac{\partial^2F}{\partial\alpha^2} (\alpha, \eta)$ is a continuous
function with respect to both $\al$ and $\eta$ we conclude that

\begin{equation*}
\lim_{n\to\infty} \frac{\partial^2F}{\partial\alpha^2}(\alpha_n,1) = 0
\Longleftrightarrow \lim_{n \to \infty} N''(\alpha_n) = 0 \ .
\end{equation*}

Since
$\alpha_n \in A$, we obtain, by Lemma~\ref{kol}(a)

\begin{equation}\lb{d1}
\frac{\partial^2F}{\partial\alpha^2} (\alpha_n,1) = \int M''(\vert f +
\alpha_n g \vert) g^2 d\mu(t) \ .
\end{equation}

If $\supp f \cap \supp g = \emptyset$ then \eqref{d1} becomes
\begin{equation*}
\frac{\partial^2F}{\partial\alpha^2} (\alpha_n,1) = \int
M''(\vert\alpha_ng\vert) g^2 d\mu(t)
\end{equation*}
and since
$\vert g\vert$ is bounded and $M''(t) \to 0$ as $t \to 0$ we conclude that
$\frac{\partial^2F}{\partial\alpha^2}(\alpha_n,1)
\operatornamewithlimits\longrightarrow 0$ as ${n\to\infty}$ and
therefore
$\lim_{n\to\infty}N''(\alpha_n) = 0$.

If $\supp f \cap \supp g \not=
\emptyset$ let $B \subset \supp f \cap \supp g$ and
$b_1, b_2, b_3, b_4>0$ be such that
$\mu (B) > 0$, $0<b_1<|f(t)|<b_2$, $0<b_3<|g(t)|<b_4$ for all $t\in B$.
Now let $n_0 \in \bbN$ be such that for all $n > n_0$ and all $t\in B$:
\begin{equation*}
M'' (|f(t)+\al_n g(t)|)) > \frac 12 M''(b_1) \ .
\end{equation*}
Then since
$M'' (\vert f(t) + \alpha_ng(t)\vert)g^2(t) \geq 0$ for a.e. $t$, we
get for $n >
n_0$:
\begin{equation*}
\int M''(\vert f + \alpha_n g\vert) g^2 d\mu(t) \geq \int_B
M''(\vert f + \alpha_n g\vert) g^2 d\mu(t)
\geq \frac 12 \mu ({B})M''(b_1)b_3^2 > 0\ .
\end{equation*}
Hence
$\lim_{n\to\infty} \frac{\partial^2F}{\partial\alpha^2}(\alpha_n,1)
\not= 0$
and therefore $\lim_{n\to\infty} N''(\alpha_n) \not= 0$.
\end{proof}

The above partial characterization of disjointness allows us to
immediately conclude that  isometries from subspaces of $L_M$
which contain enough disjointly supported bounded functions into
$L_M$ have to preserve disjointness.

\begin{theorem}\lb{isoMzero}
Assume that $M$
is an Orlicz function which satisfies condition
$\Delta_{2+}$  and such that
$M''$ is a continuous function with $M''(0)=0$.
Let $T:L_M\INTO L_M$ be an isometry. Then $T$
preserves disjointness.
\end{theorem}

\section{The case of $M''(0)=\infty$}

In this section we study Orlicz spaces $L_M$ analogous to $L_p, 1<p<2$.  In
this case we do not find any characterizations of disjointness.  Instead we
give conditions which help us to determine when support of $f$ is
contained in the support of $g$.  Our conditions do not provide a full
characterization of containment of supports but they are sufficient to
determine that isometries have to preserve the containment of supports.

\begin{prop}\lb{Minfinity}
Assume that $M$ is an Orlicz function which satisfies $\Delta_{2+}$ condition
and such that $M''$ is a continuous function on $(0,\infty)$ with
$\lim_{t\to 0}M''(t) = \infty$.  Let $f,g \in L_M$
with $\Vert f \Vert = \Vert g \Vert = 1$.

Then
\begin{itemize}
\item[(a)]  $\mu(\supp g\setminus\supp f) > 0 \Rightarrow \lim_{\alpha\to 0}
N''(\alpha) = \infty$
\item[(b)]  If $f,g$ are simple then
\begin{equation*}
\mu(\supp g \setminus \supp f) = 0 \Rightarrow \lim_{\alpha \to 0}
N''(\alpha) \not= \infty\ .
\end{equation*}
\end{itemize}
\end{prop}
\begin{proof}
Similarly as in the proof of Proposition~\ref{Mzero} we see that equation
\eqref{der} is valid i.e.
\begin{equation}\lb{der1}
\begin{split}
- \frac{\partial^2F}{\partial\alpha^2} (\alpha,N(\alpha)) &=
\frac{\partial^2F}{\partial\alpha\partial\eta}
(\alpha,N(\alpha))N'(\alpha) \\
&+\left[\frac{\partial^2F}{\partial\alpha\partial\eta}
(\alpha,N(\alpha))+\frac{\partial^2F}{\partial\eta^2}
(\alpha,N(\alpha))N'(\alpha)\right]N'(\alpha)\\
&+ \frac{\partial F}{\partial\eta}(\alpha,N(\alpha))N''(\alpha)
\end{split}
\end{equation}
Let $A$ be the set so that for
$\alpha \in A$, $N''(\al)$ exists and 
${\frac{\partial^2F}{\partial\alpha^2}(\alpha,1)}$
is given by the formula from Lemma~\ref{kol}(a).  By Lemma~\ref{kol}(c)
there exists a sequence $(\al_n) \subset A, \alpha_n \to 0$ so that
${\frac{\partial^2F}{\partial\eta\partial\alpha}
(\alpha_n, N(\alpha_n))}$ and 
${\frac{\partial^2F}{\partial\eta^2}
(\alpha_n, N(\alpha_n))}$ are
bounded for all $n \in \bbN$. $N'(\alpha)$ is a continuous function of
$\alpha$ and $\vert N'(0)\vert < \infty$ by \eqref{der0} and
Lemma~\ref{kol}(b).   Thus
${\frac{\partial^2F}{
\partial\eta\partial\alpha}(\alpha_n,N(\alpha_n))N'(\alpha_n)}$ 
and ${\frac{\partial^2F}{
\partial\eta^2}(\alpha_n,N(\alpha_n))N'(\alpha_n)}$
are bounded
for all $n \in \bbN$.  Therefore, since ${\frac{\partial F}{
\partial\eta}(\alpha_n,N(\alpha_n))}$ is bounded for all $n$ by
Lemma~\ref{kol}(b), we get from \eqref{der1} that
\begin{equation*}
\lim_{n\to\infty} \frac{\partial^2F}{\partial\alpha^2} (\alpha_n,N(\alpha_n)) =
\infty \iff \lim_{n\to\infty}N''(\alpha_n) = \infty\ .
\end{equation*}

Since ${\frac{\partial^2F}{\partial\alpha^2}(\alpha_n,\eta)}$ is 
a continuous function
with respect to $\eta$ we conclude that
\begin{equation*}
\lim_{n\to\infty} \frac{\partial^2F}{\partial\alpha^2} (\alpha_n,1) =
\infty \iff \lim_{n\to\infty}N''(\alpha_n) = \infty\ .
\end{equation*}

Since $\alpha_n \in A$, we obtain by Lemma~\ref{kol}(a)
\begin{equation*}
\frac{\partial^2F}{\partial\alpha^2} (\alpha_n,1) = \int_{\supp g}M''(\vert f
+ \alpha_ng\vert)g^2 d\mu (t)
\end{equation*}

For (a) assume that $\mu(\supp g\setminus \supp f) > 0$ and let $S \subset
\supp g \setminus \supp f$ be such that $\mu(S) > 0, \inf\{\vert
g(t)\vert : t \in S\} = s_1 > 0 $ and
$
\sup\{\vert g(t)\vert : t \in S\} = s_2 < \infty\ .
$
Denote
\begin{equation*}
m_n = \inf\{M''(|\al_n s|)\ :\ s\in[s_1,s_2]\}.
\end{equation*}

Since $\al_n\to 0$ as $n\to\infty$ and $M''(t)\to\infty$ as $t\to0$, we get
$\lim_{n\to\infty}m_n=\infty$.

Thus we have
\begin{eqnarray*}
\frac{\partial^2F}{\partial\alpha^2} (\alpha_n,1) &\geq &\int_S M''(\vert
\alpha_ng\vert)g^2 d\mu (t)\\
&\geq& \mu(S) m_n \cdot s_1^2
\operatornamewithlimits\longrightarrow_{n\to\infty} \infty .
\end{eqnarray*}
Thus $\lim_{n\to\infty}N''(\alpha_n) = \infty$ and therefore
 $\lim_{\alpha\to 0}N''(\alpha) = \infty$.

For (b) assume that $f,g $ are simple and
$\mu(\supp g \setminus \supp f) = 0$.
Set $m_f = \inf\{f(t) : t \in \supp g \}$.  Since $f$ is a simple function
$m_f > 0$ and since $g$ is bounded there exists $\alpha_0 > 0$ so that
$\vert f (t) + \alpha g(t)\vert > \frac 12 m_f$ for all $t$ and all
$\alpha$ with $\vert \alpha \vert < \alpha_0$.  Let
$\delta_0 = \sup \{M''(t) : \frac 12 m_f \leq t \leq \Vert f \Vert_\infty +
\Vert g \Vert_\infty\}$\ .
By continuity of $M''$ on $(0,\infty), \delta_0 < \infty$.  Thus for all
$n$ such that $\vert \alpha_n\vert <\alpha_0$ we have:
\begin{equation*}
\int_{\supp g} M'' (\vert f + \alpha_n g\vert)g^2 d\mu (t) \leq \delta_0
\cdot \Vert g \Vert^2_\infty < \infty
\end{equation*}
Hence
$\lim_{n\to\infty}N''(\alpha_n) \not= \infty$ and therefore
$\lim_{\alpha\to 0}N''(\alpha)\not= \infty$.
\end{proof}
\begin{cor}\lb{supp}
Assume that $M$ is an Orlicz function which satisfies $\Delta_{2+}$ condition
and such that $M''$ is a continuous function on $(0,\infty)$ with
$\lim_{t\to 0}M''(t) = \infty$.  Let $T : L_M
\INTO L_M$ be an isometry
and $f,g\in L_M$ be such that $\supp f= \supp g$ up to a set of measure
zero. Then $\supp Tf= \supp Tg$ up to a set of measure
zero.
\end{cor}

\begin{proof}
Without loss of generality we assume that
$f,g$ are simple and
$\Vert f \Vert = \Vert g \Vert = 1$.  By Proposition~\ref{Minfinity}(b) we
see that
$\lim_{\alpha\to 0}N''(\alpha) \not= \infty$.  Since $T$ is an isometry
$N(\alpha) = \Vert Tf + \alpha Tg\Vert$ and thus by
Proposition~\ref{Minfinity}(a)  $\mu(\supp Tg\setminus \supp Tf)  = 0$.

Exchanging the roles of $f$ and $g$ we symmetrically obtain $\mu(\supp
Tf\setminus\supp Tg)  = 0$.
\end{proof}

After the author presented   this paper at the Conference on Function 
Spaces at Southern Illinois University at
Edwardsville,
Abramovich and Kitover \cite{AK99} showed that an isometry $T$ between Banach function
spaces satisfies condition from the conclusion of Corollary~\ref{supp} if
and only if $T^{-1}: TX \lra X$ is disjointness preserving.  However they
also showed that in general it is possible to construct operators
satisfying the above condition but such that $TX$ does not have non-trivial
disjoint elements and, in particular, $T$ is not disjointness preserving
(cf. also \cite{AK98}).

Below we show that such a situation cannot happen in the case of Orlicz
spaces, i.e. every injective isometry preserving equality of supports does
preserve disjointness.

\begin{theorem}\lb{char}
Assume that $M$ is an Orlicz function which satisfies $\Delta_{2+}$ 
con\-di\-tion
and such that $M''$ is a continuous function on $(0, \infty)$
and $\lim_{t\to 0} M''(t) = \infty$.

Then every isometry
$T : L_M \INTO L_M$
preserves disjointness.
\end{theorem}

\begin{proof}
By Corollary~\ref{supp} it is enough to show that for any two disjoint sets $A,
B, \subset [0, 1]$ we have $\mu(\supp T\chi_A \cap \supp T\chi_B) =0$.

The proof of this fact is very short if we assume in addition that 
$M''(t)>0$ for all $t>0$ and that the complementary
function $M^*$   satisfies $\Delta_{2+}$ condition.  We present this
simple duality argument first, and then we show a longer proof which does
not require any additional assumptions.

Denote $H = \span \{\chi_A, \chi_B\}$. 
By Lemma~\ref{calvert} 
 the map $S:\span \{\chi_A^N,
\chi_B^N\}\lra\span \{(T\chi_A)^N, (T\chi_B)^N\}$ defined by
$$S\chi_A^N=(T\chi_A)^N,\ \ \ \ \  S\chi_B ^N=(T\chi_B)^N$$
is an isometry.

Notice that  for all $v$
$$(M^*)''(v) = {\frac {1}{M''\big( (M')^{-1}(v)\big) }}.$$

So $(M^*)''$ is continuous, $(M^*)''(v)>0$ for all $v>0$  and
$$\lim_{v \to 0}(M^*)''(v) = \lim_{\omega \to 0}{\frac{1}{M''(\omega)}} = 0.$$

If $M^*$ satisfies $\Delta_{2+}$-condition we obtain by
Proposition~\ref{Mzero} that $S$ preserves disjointness, since
bounded functions are clearly dense in $\span \{\chi_A^N,
\chi_B^N\}\subset L_{M^*}$. And since $\supp
h^N = \supp h$  for all $h \in L_M$, we conclude that $\mu(\supp T\chi_A
\cap \supp T\chi_B) = 0$ as desired.

If $M^*$ does not satisfy $\Delta_{2+}$ condition or if
$M''(t)\not>0$ for all $t>0$, we will apply a much longer,
more direct approach relying on the fact that 
$(TH)^{*}=\span\{(T\chi_A)^N, (T\chi_B)^N\}$ 
is a subspace of
$L_{M^*}$.

Denote $f=\chi_A,\ g=\chi_B$ and 
for any scalar $\al$ let $h=\chi_A+\al\chi_B$ and
 $Th = f + \alpha g$.  By \cite{GrzH} for almost every $t\in[0,1]$ we have:

\begin{eqnarray*}
(Th)^N(t) &=& C_\alpha M'\left(\frac{\vert(f+\alpha g)(t)\vert}{\Vert f+\alpha g
\Vert}\right) \sgn ((f+\alpha g) (t))\ ,\\
f^N(t) &=& C_f M'\left(\frac {\vert f(t)\vert}{\Vert f \Vert}\right)
\sgn (f(t)) \ ,\\
g^N(t)&=& C_g M'\left(\frac{\vert g (t)\vert}{\Vert g \Vert}\right) \sgn
(g(t))\ ,
\end{eqnarray*}
where constants $C_\alpha, C_f, C_g$ do not depend on $t$.

Thus for each $\alpha $ there exist
$\beta_1(\alpha), \beta_2(\alpha)$ so that

\begin{equation}\lb{ld}
\begin{split}
M' \left(\frac{\vert(f + \alpha g)(t)\vert}{\Vert f + \alpha
g\Vert}\right)&\sgn ((f + \alpha g)(t))=\\
  &\beta_1(\alpha) M'
\left(\frac{\vert f (t)\vert}{\Vert f \Vert} \right)
\sgn (f(t)) + \beta_2(\alpha)M' \left(\frac{\vert g(t)\vert}{\Vert g
\Vert}\right) \sgn (g(t)) .
\end{split}
\end{equation}

Since
$A \cap B = \emptyset$, by Corollary~\ref{supp},

\begin{eqnarray*}
\mu (\supp g \setminus \supp f)  & > & 0  \ \ \text{and}\\
\mu(\supp f \setminus \supp g)   & >& 0 \ .
\end{eqnarray*}
Thus for $t \in \supp f \setminus \supp g$ equation \eqref{ld}
becomes

\begin{equation}\lb{ldf}
M'\left(\frac{\vert f(t)\vert}{N(\alpha)}\right) \sgn (f(t)) =
\beta_1(\alpha) M'\left(\frac{\vert f(t)\vert}{\Vert f\Vert}\right) \sgn
(f(t)) + 0 .
\end{equation}

Our next goal is to evaluate
$\beta_1(\alpha)$.  We follow a technique similar to the one in the proof
of \cite[Theorem~6.1]{complex}.

Suppose that
$\vert f|\Big|_{\supp f\setminus\supp g}$ is not constant, that is there
exist
$t_1, t_2 \in \supp f \setminus \supp g$ so that ${\vert
f(t_1)\vert}/{\Vert f \Vert} \not= {\vert f (t_2)\vert}/{\Vert f
\Vert}$.  Denote
$x_i =  {\vert f(t_i)\vert}/{\Vert f \Vert} , \ i = 1,2$.  Then
\eqref{ldf} becomes

\begin{eqnarray*}
                   M'
                   \left(x_1 \cdot
                   \frac
                   {\Vert f \Vert}
                   {N(\alpha)}
                   \right)
                    &=&
                   \beta_1(\alpha) M'(x_1)\\
                   M'
                   \left(x_2 \cdot
                   \frac
                   {\Vert f \Vert}
                   {N(\alpha)}
                   \right)
                    &=&
                   \beta_1(\alpha) M'(x_2)
\end{eqnarray*}

Thus
\begin{equation}\lb{ldf1}
\frac{{M'\left(x_1 \cdot
 \frac{\Vert f \Vert}{N(\alpha)}
\right)}}{{M'(x_1)}}
=
\frac{M'\left(x_2 \cdot\frac{\Vert f \Vert}{N(\alpha)}\right)}
{M'(x_2)}
\end{equation}
for all $\alpha \in \bbR$.  Notice that
$N([0,\infty)) = [\Vert f \Vert, \infty)$.
Thus $({\Vert f \Vert}/{N(\alpha)}) \in (0,1]$.

Set
\begin{equation*}
u = x_2 \frac{\Vert f \Vert}{N(\alpha)},\ \ \ \beta =
\frac{M'(x_1)}{M'(x_2)},\ \ \ \xi = \frac{x_1}{x_2}\ .
\end{equation*}

Then \eqref{ldf1} becomes

\begin{equation*}
M'(\xi u) = \beta M' (u)
\end{equation*}
for all $u \in [0, x_2]$.

Thus $$M(\xi u) = \beta \xi M(u)$$
 for all $u \in [0,x_2]$ and by
\cite[Lemma~6.2]{complex}
there exists $p \geq 1$ and constants $C_1, C_2$ so that for all $u \leq
x_2$

\begin{equation*}
C_2 u^p \leq M(u) \leq C_1u^p\ .
\end{equation*}

Moreover, if $M(u) \not\equiv Cu^p$ for $u \in [0,x_2]$ then there exists
$\gamma > 0$ such that

\begin{equation}
M'(a \gamma^k) = M' (a)M' (\gamma^k)\lb{mult}
\end{equation}
for all $a > 0, k \in \bbZ$ and $\xi \in \{\gamma^k \in \bbZ\}$.

Thus if
$\vert f \vert \bigg\vert_{\supp f \setminus \supp g}$ is not constant and
if
$M(u) \not\equiv Cu^p$ for $u$ near zero then there exists $x > 0$ and
$\gamma$ satisfying \eqref{mult} so that
$ {\vert f (t)\vert}/{\Vert f \Vert} \in \{x \gamma ^k : k \in \bbZ\}$
for all
$t \in \supp f\setminus \supp g$.  Hence \eqref {ldf} becomes
 $$M'\left(x\gamma^k \frac{\Vert f \Vert}{N(\alpha)}\right) =
 \beta_1(\alpha)M' (x\gamma^k).$$

By \eqref{mult} we get
 \begin{equation*}
 M'(\gamma^k) M' \left(x \frac{\Vert f \Vert}{N(\alpha)}\right) =
\beta_1
 (\alpha) M'(x)M' (\gamma^k)\ .
 \end{equation*}

Hence
\begin{equation}
\beta_1(\alpha) =  \frac{M'\left(x\frac{\Vert
f\Vert}{N(\alpha)}\right)}{M'(x)}\ \ . \lb{beta1}
\end{equation}

Clearly, if
$\vert f \vert\bigg\vert_{\supp f \setminus \supp g} \equiv
\text{const.}
= x$ then \eqref{ldf} becomes

\begin{equation*}
M'\left(x \frac{\Vert f \Vert}{N(\alpha)}\right) = \beta_1 (\alpha) M'(x)
\end{equation*}

and \eqref{beta1} holds.

Using similar technique we show that when $M(u) \not\equiv Cu^p$ for $u$ near
zero for any $p$, then there exists $y > 0$ so that

\begin{equation}
\beta_2(\alpha) = \frac{M'\left(y \frac{\vert \alpha \vert\  \Vert g
\Vert}{N(\alpha)}\right)} {M'(y)}  \lb{beta2}\ \ .
\end{equation}

Thus for $t \in \supp f \cap \supp g,\ $ \eqref{ld} becomes:
\begin{equation}
\begin{split}
M' \left(\frac{\vert f (t) + \alpha g (t)\vert}{N(\alpha)}\right)
&\sgn (f (t) +
\alpha g (t)) = \\
 &\frac{1}{M'(x)} M'\left(x \cdot \frac{\Vert f \Vert}{N(\alpha)}\right)
\cdot M' \left(\frac{\vert f(t)\vert} {\Vert f \Vert}\right) \sgn (f(t))
 \\ &+\frac{1}{M'(y)} M'\left(y \cdot \frac{\vert \alpha| \ \Vert g
\Vert}{N(\alpha)}\right)  M' \left(\frac{\vert g (t)\vert}{\Vert g
\Vert}\right) \sgn (g(t)) \ .
 \lb{ld1}
\end{split}
\end{equation}

Fix $t \in \supp f \cap \supp g$ and let
$\alpha_0 > 0$ be small enough so that
$\sgn (f(t) + \alpha g (t)) = \sgn(f(t))$ for all
$\alpha \in (0,\alpha_0)$.  We differentiate \eqref{ld1} with respect to
$\alpha$ when
$\alpha \in (0,\alpha_0)$:

\begin{equation*}
\begin{split}
M'' &\left(\frac{\vert f (t) + \alpha g (t)\vert}{N(\alpha)}\right)
\sgn(f(t)) \cdot \frac{\sgn(f(t)) g (t)N(\alpha) - N'(\alpha)\vert f (t) +
\alpha g(t)\vert}{(N(\alpha))^2}
=\\
 &\frac{1}{M'(x)} M'\left(\frac{\vert f (t)\vert}{\Vert f \Vert}\right)
\cdot \sgn(f(t)) \cdot M'' \left(x \cdot \frac{\Vert
f\Vert}{N(\alpha)}\right) \cdot \frac{\left(- x \Vert f
\Vert\right)}{\left(N(\alpha)\right)^2}\\
&+\frac{1}{M'(y)} M' \left(\frac{\vert g (t)\vert}{\Vert g \Vert}\right)
\sgn (g(t)) \cdot M'' \left(y \cdot \frac{\alpha \Vert
g\Vert}{N(\alpha)}\right)
\cdot \frac{y \Vert g \Vert N(\alpha) - N'(\alpha)y \alpha \Vert g
\Vert}{(N(\alpha))^2}
\end{split}
\end{equation*}

When $\alpha$ approaches zero we obtain:
\begin{equation*}
\begin{split}
M''\left(\frac{\vert f (t)\vert}{\Vert f \Vert}\right) &\sgn (f(t)) \cdot
\frac {\sgn(f(t)) g (t) \Vert f \Vert - 0}{\Vert f \Vert^2} =\\
&\frac{1}{M'(x)} M'\left( \frac{\vert f (t)\vert}{\Vert f \Vert}\right)
\sgn (f(t)) M''(x) \cdot \left(\frac{-x\Vert f \Vert}{\Vert f
\Vert^2}\right)\\
&+ \frac {1}{M'(y)} M'\left(\frac{\vert g (t)\vert}{\Vert g \Vert}\right)
\sgn(g(t)) \cdot \frac {y \Vert g \Vert\cdot \Vert f \Vert - 0}{\Vert f \Vert^2}
\cdot \lim_{\alpha \to 0} M''\left(\alpha \frac{y \Vert g \Vert}{N(\alpha)}
\right)
\end{split}
\end{equation*}
which is a contradiction since

$$
\lim_{\alpha \to 0} M''\left(\alpha \cdot \frac{y\Vert g
\Vert}{N(\alpha)}\right) = \infty
$$
and all other quantities in the above equation are finite and nonzero.

Hence $\mu (\supp f \cap \supp g) = 0$.
\end{proof}

\section{Final remarks}
In this section we summarize the results that we obtained:

\begin{cor}\lb{final}
Suppose that $M$ is an Orlicz function satisfying condition $\Delta_{2+}$ and
such that $M''$ is continuous and either $\lim_{t \to 0}M''(t)=\infty$ or
$M''(0)=0$ and $M''(t)>0$ for all $t>0$.
Suppose further that $T:L_M \INTO L_M$ is an isometry.
Then there exist a Borel map $\sigma:[0,1]\lra[0,1]$ and a function $a:[0,1]
\lra \bbR$ (or, if $L_M$ is complex, $a:[0,1] \lra \bbC$) so that for every
$f \in L_M$ and almost every $t \in [0,1]$:

\begin{equation}\lb{star}
Tf(t) = a(t)f (\sigma(t))
\end{equation}

Moreover $\vert a(t) \vert = 1$ a.e. unless there exist constants $C_1,
C_2, t_0 > 0$ and $p, 1<p<\infty$,

so that for all $t \leq t_0$:
$$C_1 t^p \leq M(t) \leq C_2 t^p.$$

If such constants exist, but $(M(t)/t^p) \not\equiv const.$ on any
interval containing $0$ then there exist $A, \gamma > 0$ so that for a.e.
$t$,
$\vert a(t) \vert = A \cdot \gamma^{k(t)}$, where $k(t)\in \bbZ$.
\end{cor}

\begin{proof}
It follows from Theorems~\ref{isoMzero} and \ref{char} that $T$ preserves
disjointness.  Abramovich \cite{Ab83} proved that this implies that $T$ is
a weighted composition operator i.e. $T$ has form \eqref{star}.

To prove the "moreover part" we will use the nonatomic version of
\cite[Theorem~6.1]{complex} which we remind below:

\begin{theorem}\lb{complex}
Let $M$ be a  continuous Orlicz function and let $f,g \in L_M$ be disjoint
elements such that $f, g \not\equiv 0$ and 
$\{h^N\ :\ h\in\span\{f,g\}\}=\span\{f^N,g^N\}\subset L_M^*$.  
Then one of three possibilities holds:

\begin{itemize}
\item[(1)] there exists a scalar $k_f$ so that $\vert f(t)\vert = k_f$ for
almost all $t \in supp f$; or
\item[(2)] there exists $p, 1 \leq p < \infty$, $C \geq 0$ and $t_0 > 0$ so
that $M(t) = Ct^p$ for all $t \leq t_0$; or
\item[(3)] there exist $p, 1 \leq p< \infty$ and constants $C_1, C_2,
t_0,\gamma, k_f \geq 0$ so that $ C_2 t_p \leq M(t)\leq C_1 t^p$ for all $t
\leq t_0$ amd such that for almost all $t \in supp f$ there exists $k(t)\in
\bbZ$ with $$\vert f(t)\vert = k_f{\cdot}\gamma^{k(t)}$$
\end{itemize}

\end{theorem}

\begin{rem}
\cite[Theorem~6.1]{complex} is stated and proven for sequence spaces
$\ell_M$, but the nonatomic version requires only very minor routine
adjustments, so we leave them to the interested reader.
\end{rem}

Now we recall that by Lemma~\ref{calvert}, for any disjoint subsets 
$A,B$ of $[0,1]$ we have
$\{(Th)^N\ :\ h\in \span\{\chi_A,\chi_B\}\}
=\span \{(T\chi_A)^N, (T\chi_B)^N\}$.
  But we know that $T\chi_A, T\chi_B$ are
disjointly supported so the "moreover part" of Corollary~\ref{final}
follows directly from Theorem~\ref{complex} applied to $f = T\chi_A, g =
T\chi_B$.
\end{proof}

\begin{rem}
 The statement of Corollary~\ref{final} leaves open the case when
$M''$ is continuous and $M''(0)=c$ for some $c,$ $ 0<c<\infty$.  The function $M_2(t)=t^2$ belongs to
this case and of course $L_{M_2}$ has non-disjointness preserving injective
and surjective isometries.

If $M(t) \not\equiv t^2$ on $[0,1]$ then it is known that all surjective
isometries are disjointness preserving; however our differential technique
does not seem to provide enough information about injective isometries in
this case.  We feel that the hardest case would be to distinguish behavior
in $L_2$ from $L_M$ where $M(t)= t^2$ for all $t \leq a < 1$ but $a$ is
close to 1.
\end{rem}

\begin{rem}
Our results deal with injective isometries where domain is entire
$L_M$(see also the remark before Theorem~\ref{isoMzero}).  It would be
interesting to determine if isometries from subspaces of $L_M$ into $L_M$
have to be disjointness preserving, as it is the case in $L_p, p \neq 2$
(cf. \cite{Kol91}); (note that when $M(t)= t^2$ for all $t \leq a < 1$,
where $a$ is large enough, then $L_M$ contains an isometric copy of
$\ell_2^2$ \cite[Example~3]{complex}, so clearly injective isometries from
the subspace of $L_M$ do not have to preserve disjointness in this case).
\end{rem}

%\bibliographystyle{plain}
%\bibliography{tref}

\end{document}